\documentclass[12pt]{article}

\usepackage[english] {babel}
\usepackage{amsmath,amssymb,amscd}
\usepackage{amsthm}
\usepackage{verbatim}
\usepackage[utf8]{inputenc}
\usepackage{graphicx}
\usepackage{fancybox}
\usepackage{subfigure}
\usepackage[usenames,dvipsnames]{color}
\usepackage{sidecap}
\usepackage{caption}
\usepackage{ifthen}
\usepackage{textcomp}
\usepackage[all]{xy}
\usepackage{feynmf}
\usepackage{anysize}

\newcommand{\A}{{\mathbb A}}
\newcommand{\N}{{\mathbb N}}

\newcommand{\Z}{{\mathbb Z}}
\newcommand{\R}{{\mathbb R}}
\newcommand{\C}{{\mathbb C}}
\newcommand{\Q}{{\mathbb Q}}
\newcommand{\I}{{\mathcal{O}}}

\theoremstyle{plain}
\newtheorem{teo}{Theorem}[section]
\newtheorem{lema}[teo]{Lemma}
\newtheorem{cor}[teo]{Corollary}
\newtheorem{prop}[teo]{Proposition}

\theoremstyle{definition}
\newtheorem{defi}{Definition}[section]

\theoremstyle{remark}
\newtheorem{obs}{Remark}[section]

\def\noi{\noindent}

\begin{document}

\nocite{*}

\title{Adelic solenoid I: Structure and topology}

\author{Juan M. Burgos, Alberto Verjovsky \thanks{Grant support listed here.}}

\maketitle

\noi {\bf 2010 Mathematics Subject Classification.} Primary: 11R56, 32L05, 55M25, 57R30.
\noi Secondary: 22Bxx.


\begin{abstract} 
\noi Topologically the adelic Riemann sphere is the suspension of the adelic solenoid and because of this relation, here we study the adelic solenoid by studying the adelic Riemann sphere topology. The main result is the Birkhoff-Grothendieck Theorem: A holomorphic vector bundle splits as a sum of holomorphic line bundles whose Chern character is now a rational number. As a consequence, the Picard group is isomorphic to the additive group of rational numbers and the $K$--ring has new elements that factor the tautological class.
\end{abstract}

%\noindent { \bf Key Words:} Adelic solenoid, adelic quasi-conformal map, adelic quasi-symmetric homemorphism
%and Teichm\"uller space.\\
%{\bf 2010 Mathematics Subject Classification.} Primary: 11R56, 32G05, 32G15, 57R30.
%Secondary: 22Bxx.

%\tableofcontents

\section{Introduction}
Since its creation by Claude Chevalley and Andr\'e Weil the ring of ad\`eles of the rationals $\A_\Q=\R\times\underset{p}\prod'\,\Q_p$  has played a fundamental role in number theory, for instance in {\em class field theory} \cite{RV}, \cite{Ta} and the Langlands program. The canonical diagonal inclusion $i:\Q\to\A_\Q$ embeds $\Q$ as a discrete cocompact subgroup of 
$\A_\Q$ which we identify with $\Q$. The quotient $\A_\Q/\Q$ is the {\em ad\`ele class group} with its additive structure is a compact abelian group and its Pontryagin dual is the additive group of the rationals $(\Q,+)$ with the discrete topology. 
There is another description of $\A_\Q/\Q$ as the inverse limit of all finite coverings $p_n(z)={z} ^n$  ($z\in{S}^1,n\in\Z$) of the circle $S^1$. This is a one dimensional solenoidal compact abelian group in the sense of Pontryagin.  It is a sort of ``diffuse circle" 
(a {\em lamination}, a {\em current} or  a {\em foliated cycle} in the sense of Sullivan \cite{Su2}) and in fact we denote this group as 
$$
S^1_\Q=\Q^{\vee}=\text{Pontryagin\, dual of} \,\, \Q
$$
to convey the idea that it is a generalization of a circle. This solenoid is a lamination with dense leaves which are densely immersed copies of the real line. The solenoid $S^1_\Q$ is also called the \textsf{algebraic universal covering space} of the circle $S^1$. The Grothendieck--Galois group of the covering is $\hat{\Z}$, the \textsf{\'etale fundamental group} of $S^1_\Q$ \cite{SGA1}.

Concerning dynamical systems theory, in \cite{Universalities} (see also \cite{Su2}) D. Sullivan  studies the linking between universalities of Milnor-Thurston, Feigenbaum's (quantitative) and Ahlfors-Bers. As he points out in his second example, the dynamical suspension of the square map from the diadic solenoid $S^{1}_{2}$ to itself is the basic solenoidal surface required in the dynamical theory of Feigenbaum's Universality \cite{Feigenbaum}. This object can be seen as a quotient of $\Delta^{*}_\Q$, the inverse limit of the punctured disk by the inverse system of coverings considered before. This is a two dimensional solenoid with hyperbolic leaves and so is $L$. This work was continued, for instance, in the use of 3-dimensional hyperbolic laminations by Misha Lyubich and Yair Minsky in \cite{LM}.

If instead of the circle we consider the same inverse system of coverings on the Riemann sphere, its inverse limit is the \textsf{adelic Riemann sphere} and we denote it by $\C P^{1}_\Q$. Now the finite coverings ramify on $0$ and $\infty$ giving two cusps in the inverse limit. Topologically, it is the suspension of the adelic solenoid and because of this relation, here we study the adelic solenoid by studying the adelic Riemann sphere topology. We prove the analog of the Birkhoff-Grothendieck Theorem \cite{Birkhoff}, \cite{Grothendieck}: A holomorphic vector bundle splits as a sum of holomorphic line bundles whose Chern character is now a rational number. In particular, the Picard group is isomorphic to the additive group of rational numbers instead of the integers as in the usual case. The $K$--ring has also new elements that factor the tautological class $x$:
$$K\left(\C P^{1}_\Q\right)\cong \frac{\Z[x^{q}\ /\ q\in\Q]}{\left\langle\ (x^{q}-1)(x^{p}-1)\ /\ q,p\in\Q \right\rangle}$$
Note that, as it was expected, the $K$--ring of the Riemann sphere embeds in this ring:
$$\frac{\Z[x]}{\left\langle\ (x-1)^{2} \right\rangle}\cong K(\C P^{1}) \hookrightarrow K\left(\C P^{1}_\Q\right)$$
where $x$ is the image of the class represented by the tautological complex line bundle.

For the proof of these results, first we must state some preliminaries about the adelic solenoid algebraic structure and limit periodic maps. We also develop the notion of degree for a continuous map from the adelic solenoid to itself. It is a homotopy invariant taking values in the rational numbers and extends the usual degree in circles.

This is the first part of a series of two papers. The next part generalizes Ahlfors-Bers theory to the solenoidal context.

\section{Adelic solenoid}\label{adelic_solenoid}

In what follows we will identify the group $U(1)$ with the unit circle $S^{1}= \{z\in\C:\, |z|=1\}$ and the finite cyclic group $\Z/n\Z$ with the group of $n^{th}$ roots of unity in $S^1$. 

By covering space theory, for any integer $n\neq 1$, it is defined the unbranched covering space of degree $n$, $p_n:S^1 \to S^1$ given by $z\longmapsto {z^n}$. If $n,m\in \Z^+$ and $n$ divides $m$, then there exists a covering map $p_{n,m}:S^1\to S^1$ such that $p_n \circ {p_{n,m}} = p_m$ where 
$p_{n,m}(z)=z^{m/n}$. 
We also denote with the same letters the restriction of $p_n$ and $p_{n,m}$ to the $n^{th}$ roots of unity. In particular we have the relation:
$$p_{n,m} \circ{p_{m,l}}=p_{n,l}$$
 This determines a projective system of covering spaces $\{S^1,p_{n,m}\}_{n,m \geq 1, n|m}$ whose projective limit is the \textsf{universal one--dimensional solenoid} or \textsf{adelic solenoid}
\[ S^1_\Q:=\lim_{\overset{p_n}{\longleftarrow} }S^1. \] 
Thus, as a set, $S^1_\Q$ consists of sequences $\left( z_n\right)_{n\in\N,\, z\in{S^1}}$ which are compatible with $p_n$ i.e. $p_{n,m}(z_m)=z_n$ if $n$ divides $m$.

The canonical projections of the inverse limit are the functions  $ S^1_\Q\overset{\pi_n}\to S^1$ defined by $\pi_n\left(\left(z_j\right)_{j\in\N}\right)=z_n$ and they define the solenoid topology as the initial topology of the family. The solenoid is an abelian topological group and each $\pi_n$ is an epimorphism.  In particular each $\pi_n$ is a character which determines a locally trivial $\hat{\Z}$--bundle structure where the group
$$\hat{\Z}:=\lim_{\overset{p_n}{\longleftarrow} }\Z/m\Z$$ 
is the profinite completion of $\Z$, which is a compact, perfect and totally disconnected abelian topological group homeomorphic to the Cantor set. 
Being $\hat{\Z}$ the  profinite completion of $\Z$, it admits a canonical inclusion of $\Z\subset\hat\Z$ whose image is dense.
We have an inclusion $\hat\Z\overset{\phi}\to{S^1_\Q}$ and a short exact sequence
$0\to{\hat\Z}\overset{\phi}\rightarrow S^1_\Q\overset{\pi_1} \rightarrow S^1\to1$.

The solenoid $S^1_\Q$ can also be realized as the orbit space of the $\Q$--bundle structure $\Q \hookrightarrow \mathbb{A} \to \A/\Q$, where $\A$ is the ad\`ele group of the rational numbers which is a locally compact Abelian group, $\Q$ is a discrete subgroup of $\A$ and $\A/\Q \cong S^1_\Q$ is a compact Abelian group (see \cite{RV}). From this perspective, $\A/\Q$ can be seen as a projective limit whose $n$--th component corresponds to the unique covering of degree $n\geq 1$ of $S^1_\Q$.

By considering the properly discontinuously diagonal free action of $\Z$ on $\hat{\Z}\times\R$  given by
\[ n\cdot(x,t)=(x+n,t-n), \quad (n\in \Z, \, x\in\hat\Z, \,t\in\R)\]
the solenoid $S^1_\Q$ is identified with the orbit space $\hat{\Z}\times_{\Z} \R$. Here, $\Z$ is acting on $\R$ by covering transformations and on $\hat{\Z}$ by translations. The path--connected component of the identity element $1\in S^1_\Q$ will be called the \textsf{baseleaf} \cite{Odden} and it is a densely immersed copy of $\R$.
\smallskip

Hence $S^1_\Q$ is a compact, connected, abelian topological group and also a one-dimensional lamination where each ``leaf" is a simply connected one-dimensional manifold, homeomorphic to the universal covering space $\R$ of $S^1$, and a typical ``transversal" is isomorphic to the Cantor group $\hat{\Z}$. The solenoid $S^1_\Q$ also has a leafwise $\mathrm{C}^\infty$ Riemannian metric (i.e., $\mathrm{C}^\infty$ along the leaves) which renders each leaf isometric to the real line with its standard metric $dx$. So, it makes sense to speak of a rigid translation along the leaves. The leaves also have a natural order equivalent to the order of the real line hence also an orientation.

Summarizing the above discussion we have the commutative diagram: 
\begin{equation}\label{diagram_I}
\xymatrix{
    S^{1}_{\Q}= \lim S^{1}  \ar[r] &	\ldots S^{1} \ar[r]^{p_{m,n}}		& 	S^{1} \ar[r] 		& 	\ldots S^{1}   \\
	\hat{\Z}= \lim \Z/n\Z \ar@{^{(}->}[u]^{\phi} \ar[r]  & 	\ldots \Z/n\Z \ar@{^{(}->}[u]_ {l\,\mapsto{e^{2\pi il/n}}} \ar[r]^{p_{m,n}} 	&	\Z/m\Z \ar@{^{(}->}[u]_{l\,\mapsto{e^{2\pi il/m}}} \ar[r]	&	 \ldots \{0\} \ar@{^{(}->}[u]_{0\,\mapsto1}}
\end{equation}

\noi where $\hat{\Z}$ is the adelic profinite completion of the integers and the image of the group monomorphism $\phi:(\hat{\Z},+)\rightarrow (S^{1}_{\Q},\cdot)$ is the \textsf{principal fiber}. We notice that $\pi_{n}(x)= \pi_{n}(y)$ implies $\pi_{n}(y^{-1}x)=1$ and therefore $y^{-1}x= \phi(a)$ where $a\in n\hat{\Z}$ for some $n\in\Z\subset\hat{\Z}$.

\begin{lema}\label{Exact_seq_I}
The following is a short exact sequence of topological abelian groups:

$$\xymatrix{ 0 \ar[r] & \hat{\Z} \ar[r]^{\phi} & S^{1}_{\Q} \ar[r]^{\pi_{1}} & S^{1} \ar[r] & 1}$$

\noi and we have the commutative diagram:

$$\xymatrix{ 	0 \ar[r] & \hat{\Z} \ar[r]^{\phi} & S^{1}_{\Q} \ar[r]^{\pi_{1}} & S^{1} \ar[r] & 1 \\
			0 \ar[r] & n\hat{\Z} \ar@{^{(}->}[u] \ar[r]^{\phi} & S^{1}_{\Q} \ar[r]^{\pi_{n}} \ar[u]_{=} & S^{1} \ar[r] \ar[u]_{p_n} & 1}$$

\end{lema}
\begin{proof}
\begin{itemize}

\item By definition the following diagram commutes:

$$\xymatrix{
    S^{1}_{\Q}  \ar[r]^{\pi_{1}} &	S^{1}   \\
	\hat{\Z} \ar@{^{(}->}[u]^{\phi} \ar[r]  & 	 \{0\} \ar@{^{(}->}[u]_{0\,\mapsto1}}$$

In particular $\pi_{1}\circ \phi= 1$ and $Im(\phi)\subset  \ker(\pi_{1})$. Suppose that $\pi_{1}(x)=1$. Then
$$x=(\ldots, a_{n},\ldots, a_{m},\ldots 1)= (\ldots, e^{2\pi i b_{n}/n},\ldots, e^{2\pi i b_{m}/m},\ldots 1)= \phi(y)$$
such that $y=(\ldots, b_{n},\ldots, b_{m},\ldots 0)$. We have proved that $\ker(\pi_{1})\subset Im(\phi)$. Because $\pi_{1}$ is an 
epimorphism and $\phi$ is a monomorphism we have the first item.

\item For the second item, the second exact sequence follows exactly from the same arguments as the first. Because $\pi_{1}= z^{n}\circ \pi_{n}$, we have the right commutative square. The left square is trivial (diagram chasing).
\end{itemize}
\end{proof}

We define the \textsf{baseleaf} as the image of the monomorphism $\nu:\R \rightarrow S^{1}_{\Q}$ defined as follows:
\begin{equation}\label{diagram_II}
\xymatrix{
    S^{1}_{\Q}= \lim S^{1}  \ar[r] &	\ldots S^{1} \ar[r]^{p_{m,n}} 		& 	S^{1} \ar[r] 		& 	\ldots S^{1}   \\
	\R \ar@{^{(}->}[u]_{\nu} \ar[r]^{=}  & 	\ldots \R \ar[u]_{t\,\mapsto{e^{it/n}}} \ar[r]^{=} 	&	\R \ar[u]_{t\,\mapsto{e^{it/m}}} \ar[r]^{=}	&	 \ldots \R \ar[u]_{t\,\mapsto{e^{it}}}}
\end{equation}
	
In particular, the immersion $\nu$ is a group morphism and comparing the diagrams \eqref{diagram_I} and \eqref{diagram_II}, we have $\nu(2\pi n)= \phi(n)$ for every integer $n$. Define: $$\exp:\ \hat{\Z}\times\R\rightarrow S^{1}_{\Q}$$ such that $\exp(a,\theta)= \phi(a).\nu(\theta)$.

\begin{lema}\label{Exact_seq_II}
The following is a short exact sequence of topological abelian groups:

$$\xymatrix{ 0 \ar[r] & \Z \ar[r]^{\iota} & \hat{\Z}\times \R \ar[r]^{\exp} & S^{1}_{\Q} \ar[r] & 1}$$
such that $\iota(a)= (a,-2\pi a)$.
\end{lema}
\begin{proof}
\textit{The exponential is an epimorphism:} Consider $x\in S^{1}_{\Q}$ and $a\in \R$ such that $e^{ia}= \pi_{1}(x)$. Because $\pi_{1}\circ \nu= e^{i\theta}$ we have that $\pi_{1}(\nu(a))= e^{ia} = \pi_{1}(x)$; i.e. $\pi_{1}(\nu(a)^{-1}x)=1$. By Lemma \ref{Exact_seq_I} there is an adelic integer $b\in \hat{\Z}$ such that $\phi(b)= \nu(a)^{-1}x$; i.e. $x= \phi(b)\nu(a)= \exp(b,a)$.

\textit{The exponential kernel:} Suppose that $\exp(a,\theta)= \phi(a)\nu(\theta)=1$. Then $\phi(a)= \nu(-\theta)$ and composing with $\pi_{1}$ we have $1= e^{-i\theta}$ and $\theta= 2\pi k$ for some integer $k$. Then $$1=\exp(a,2\pi k)= \phi(a)\nu(2\pi k)= \phi(a)\phi(k)= \phi(a+k)$$
Because $\phi$ is monomorphism we have that $a+k=0$. We conclude that $a$ is an integer and $\theta= -2\pi a$
\end{proof}

\begin{cor}

\begin{enumerate}
\item $\pi_{1}:S^{1}_{\Q}\rightarrow S^{1}$ is a fiber bundle with fiber isomorphic to $\hat{\Z}$ and monodromy the shift $T(x)=x+1$.
\item $\exp$ is a local homeomorphism.
\item Restricted to a leaf, $\pi_{1}$ is a local homeomorphism.
\item $S^{1}_{\Q}$ is the dynamical suspension of the shift $T(x)=x+1$.
\item $S^{1}_{\Q}$ is foliated by dense $\R$-leaves.
\end{enumerate}
\end{cor}
\begin{proof}
\begin{enumerate}
\item If $diam(U)<2\pi$ then $U$ is a trivializing neighborhood of $S^{1}$.
\item $\Z$ acts as translations by $\iota(\Z)$ and because $\iota(\Z)$ is discrete in $\hat{\Z}\times \R$ then $\Z$ acts proper and discontinuously. We conclude that $\exp$ is a local homeomorphism. 
\item By definition $\pi_{1}$ is an open continuous epimorphism. Restricted to a leaf and a trivializing neighborhood $\pi_{1}$ is one to one.
\item $(x,2\pi) +\iota(1)= (x+1, 0)$ so $(x,2\pi)\sim (x+1,0)$.
\item The foliation $\hat{\Z}\times \R$ is invariant under translations by $\iota(a)$ for every integer $a$ hence it induces a foliation in the solenoid. $\Z$ is dense in its profinite completion $\hat{\Z}$ and so is every coset of $\hat{\Z}/\Z$. By the preceding item, we have that every $\R$-leaf is dense in the solenoid.
\end{enumerate}
\end{proof}

\section{Limit periodic functions}

\begin{defi}
Consider a metric space $(X,d)$ and a function $f:\R\rightarrow X$. $f$ is limit periodic if for every $\varepsilon>0$ there is a natural number $N$ such that $d\left(f(x+2\pi n), f(x)\right)<\varepsilon$ for every $x\in \R$ and every $n\in N\Z$.
\end{defi}

Do not confuse this concept with the more general one of almost periodic functions based on relatively dense sets \cite{Bohr}. An interesting discussion relating limit periodic functions, solenoids and adding machines can be found in \cite{Bell}. Usually the concept is defined for real or complex valued functions. However, we will need the generality of considering a metric space.

\begin{lema}
Consider a metric space $(X,d)$ and a function $f:\R\rightarrow X$. There is a unique function $\hat{f}$ continuous on the fibers such that the family of its restrictions to the fibers is uniformly equicontinuous and $f= \hat{f}\circ\nu$ iff $f$ is limit periodic.%  the following diagram commutes:
%$$\xymatrix{ S^{1}_{\Q} \ar[r]^{\hat{f}}  & X \\
%\R \ar@{^{(}->}[u]_{\nu} \ar[ur]_{f} & }$$
%Conversely, if there is a function $\hat{f}$ with the properties above, then $f$ is limit periodic.
\end{lema}
\begin{proof}
Suppose that $f$ is limit periodic. For every real number $x$, define the map $g_x:\Z\rightarrow X$ such that $g_x(n)= f(x+2\pi n)$. We will prove that $g_x$ is uniformly continuous respect to the profinite topology on $\Z$ hence it admits a unique extension to $\hat{\Z}$.

Consider $\varepsilon>0$. There is a natural number $N$ such that $d\left(f(y+2\pi n), f(y)\right)<\varepsilon$ for every $y\in \R$ and every $n\in N\Z$. In particular, for every real number $x$ we have:
$$d\left( g_x(n), g_x(m)\right)= d\left(f(x+2\pi n), f(x+2\pi n +2\pi (m-n))\right)<\varepsilon$$
for every pair of integers such that $m-n\in N\Z$; i.e. The family $g_x$ is uniformly equicontinuous respect to the profinite topology. Consider its unique continuous extension $\hat{g}_x$ to $\hat{Z}$ and define $F:\hat{\Z}\times \R\rightarrow X$ such that $F(z,x)= \hat{g}_x(z)$. Uniqueness of the extension implies that the relation $g_x(n+m)= g_{x+2\pi m}(n)$ extends to $\tilde{g}_x(z+m)= \tilde{g}_{x+2\pi m}(z)$ and we have $F(z+m, x)=F(z, x+2\pi m)$. By Lemma \ref{Exact_seq_II}, $F$ defines a unique function $\tilde{f}$ on the solenoid and by construction it is continuous on the fibers.

Conversely, given $\varepsilon>0$ there is a natural number $N$ such that:
$$d\left(\hat{f}\circ\exp(n,x), \hat{f}\circ\exp(m,x)\right)<\varepsilon$$
for every real number $x$ and every pair of integers such that $m-n\in N\Z$. Then:
$$d\left(f(x+2\pi n), f(x)\right)= d\left(\hat{f}\circ \nu(x+2\pi n), \hat{f}\circ\nu(x)\right)=d\left(\hat{f}\circ\exp(n,x), \hat{f}\circ\exp(0,x)\right)<\varepsilon$$
for every real number $x$ and integer $n\in N\Z$; i.e. $f$ is limit periodic.
\end{proof}

\begin{cor}
A continuous limit periodic function is uniformly continuous.
\end{cor}
\begin{proof}
In the previous Lemma, if $f$ is continuous then, because $\nu(\R)$ is dense in the solenoid, $\hat{f}$ is also continuous. The solenoid is compact hence $\hat{f}$ is uniformly continuous and so is $f$.
\end{proof}

\begin{lema}\label{Aprox}
Consider a metric space $(X,d)$ and a function $f$ from the solenoid to $X$ such that its restrictions to the fibers is a uniformly equicontinuous family. For every $\varepsilon>0$, there is a function $f'$ such that $d(f,f')<\varepsilon$ and it factors through some canonical projection $\pi_n$.
Moreover, if $f$ is continuous and $(X,d)$ is a locally euclidean, we can choose $f'$ to be continuous.
\end{lema}
\begin{proof}
There is a natural number $N$ such that:
$$f\left(\exp(z+N\hat{\Z}, x)\right)\subset B_{\varepsilon}\left(f(\exp(z,x))\right)$$
for every $z\in \hat{\Z}$ and $x\in \R$. Consider the lifting $s$ of $\pi_N$ such that $s(1)=e$, $s$ is continuous on $S^{1}-\{1\}$ and is continuous from one of the two directions at $1$. Define $f':=f\circ s\circ \pi_N$.

If $f$ is continuous, taking $N$ big enough we may suppose that $\exp(N\hat{\Z}, 0)\subset U$ where $U$ is a local euclidean chart at $e$. Consider a continuous function $h:X\rightarrow X$ such that $h\left(f(\exp(N\hat{\Z}, 0))\right)= f(e)$ and is the identity outside $U$. Now, take as $f'$ the function $h\circ f\circ s\circ \pi_N$. The function $h$ corrects the discontinuity at $1$ preserving the condition $d(f,f')<\varepsilon$.
\end{proof}

\begin{cor}
\begin{enumerate}
\item Every limit periodic function is the uniform limit of periodic functions.
\item If $(X,d)$ is locally euclidean, then every continuous limit periodic function is the uniform limit of continuous periodic functions.
\end{enumerate}
\end{cor}
\begin{proof}
Follows from direct application of the Lemmas above.
\end{proof}

\section{Degree}

The baseleaf map $\nu:\R\rightarrow S^{1}_\Q$ induce a coarser topology on the real line such that $\nu$ becomes an embedding. This will be called the \textsf{leaf topology} and it is generated by the following basic open sets:
$$U(x, \varepsilon, N):= \bigcup_{n\in\Z} I(x+2\pi Nn, \varepsilon)$$
where $I(x,\varepsilon)$ is the interval with center $x$ and radius $\varepsilon$. The real line with the leaf topology will be denoted by $\R_{L}$. Every continuous map from $\R_{L}$ to itself is also continuous from $\R$ to itself. However, the converse is not true.

A continuous map $f$ from the solenoid to itself can be assumed to be baseleaf preserving just by multiplying it by $f(1)^{-1}$.

\begin{lema}\label{LimitPer}
Consider a continuous baseleaf preserving function $f$ from the adelic solenoid to itself. There is a unique rational number $q$ and a unique continuous limit periodic function $h$ such that $f_{0}(x)= qx + h(x)$ where $f_{0}$ is defined by the following commutative diagram:
%\begin{equation}\label{diagram}
$$\xymatrix{S^{1}_\Q \ar[r]^{f} & S^{1}_\Q \\
			\R \ar[r]^{f_{0}} \ar@{^{(}->}[u]_{\nu} & \R \ar@{^{(}->}[u]_{\nu}}$$
%\end{equation}
\end{lema}
\begin{proof}
Respect to the leaf topology, the baseleaf is an embedding and because the solenoid is compact we have that $f_0$ is uniformly continuous respect to the leaf topology; i.e. For every $\varepsilon >0$ and natural number $\lambda$ there is a real number $\delta>0$ and a natural number $N$ such that:
\begin{equation}\label{BL_cont}
f_{0}\left( U(x,\delta, N)\right)\subset U\left(f_0(x), \varepsilon, \lambda\right)
\end{equation}
for every $x\in \R$.

Consider $\varepsilon<\pi/2$ and define $g_{m}:\R\rightarrow \R$ such that $g_{m}(x)= f_{0}(x+ 2\pi Nm )- f_{0}(x)$ for every integer $m$. We will prove that there is an integer $k$ such that $g_{m}(\R)\subset I(2\pi km,\varepsilon)$ for every integer $m$. We will prove it in the following steps:
\begin{itemize}
\item \textit{Base case:} Because of $\eqref{BL_cont}$ we have that $g_1(\R)\subset U(0,\varepsilon, \lambda)$. Because $U(0,\varepsilon, \lambda)$ is a disjoint union of open intervals, $g_1$ is continuous and $\R$ is connected, the image of $g_1$ must be contained in one of these intervals; i.e. There is an integer $k$ such that $g_{1}(\R)\subset I(2\pi k,\varepsilon)$ and $\lambda|k$.

\item \textit{Induction step:} Suppose that $g_{m}(\R)\subset I(2\pi km,\varepsilon)$ for every natural $m\leq M$. Because $g_{M+1}(x)= g_{M}(x+2\pi N)+g_{1}(x)$ and the inductive hypothesis, we have that $g_{M+1}(\R)\subset I(2\pi k(M+1) ,\pi)$. By equation \eqref{BL_cont} we have $g_{m}(\R)\subset U(0,\varepsilon, \lambda)$ for every integer $m$. Then,
$$g_{M+1}(\R)\subset I(2\pi k(M+1) ,\pi)\cap U(0,\varepsilon, \lambda)= I(2\pi k(M+1) ,\varepsilon)$$
\item \textit{Trivial case:} $g_{0}(\R)= \{0\}\subset I(0,\varepsilon)$.
\item \textit{Negative integers:} $g_{-m}(\R)= -g_{m}(\R)\subset -I(2\pi km,\varepsilon)= I(2\pi k(-m),\varepsilon)$ for every natural $m$.
\end{itemize}

We have proved a stronger version of equation \eqref{BL_cont}: For every $\varepsilon >0$ and natural number $\lambda$ such that $\varepsilon<\pi/2$, there is a real number $\delta>0$, a natural number $N$ and an integer $k$ such that:
\begin{equation}\label{cont_baseleaf2}
f_{0}\left( I( x+ 2\pi Nm, \delta)\right)\subset I\left( f_{0}(x)+ 2\pi km,\varepsilon\right)
\end{equation}
for every $x\in \R$ and every integer $m$.

Let's see that the quotient $k/N$ is independent of the $\varepsilon$ and $\lambda$ chosen. Consider another $0<\varepsilon'<\pi/2$ and $\lambda'$. There is a real number $\delta'>0$ such that $\delta'<\delta$, a natural number $N'$ and an integer $k'$ such that:
\begin{equation}\label{cont_baseleaf3}
f_{0}\left( I( x+ 2\pi N'm', \delta')\right)\subset I\left( f_{0}(x)+ 2\pi k'm',\varepsilon'\right)
\end{equation}
for every $x\in \R$ and every integer $m'$. Choose $m$ and $m'$ such that $N'm'=Nm$. Then,
$$\emptyset \neq f_{0}\left( I(2\pi N'm', \delta')\right)\subset I\left( f_{0}(0)+ 2\pi km,\varepsilon\right)\cap  I\left( f_{0}(0)+ 2\pi k'm',\varepsilon'\right)$$ and because $\varepsilon,\varepsilon'<\pi/2$ we have that $k.m=k'.m'$ hence $k/N= k'/N'$. Denote this $\varepsilon,\lambda$--independent rational by $q$. We claim that:
$$f_{0}(x)= q x\ +\ h(x)$$
where $h$ is a continuous limit periodic function. In effect, because $f_{0}$ is continuous we have that $h$ is continuous. It rest to show that it is limit periodic. Because we proved that the rational $q$ was $\varepsilon,\lambda$-independent, equation \eqref{cont_baseleaf2} implies the following: For every $\varepsilon>0$ there is a natural number $N$ such that:
$$h(x+ 2\pi Nm)-h(x)= f_{0}(x+ 2\pi Nm)- f_{0}(x)- 2\pi q\ Nm \in I(0,\varepsilon)$$
for every $x\in \R$ and every integer $m$. This proves the claim.

Moreover, this decomposition is unique for a linear limit periodic function must be zero.
\end{proof}

\begin{defi}
The rational number in Lemma \ref{LimitPer} will be called the degree of $f$ and will be denoted $\deg(f)$.
\end{defi}

\begin{prop}\label{Hom_Invariant}
The degree is a homotopic invariant.
\end{prop}
\begin{proof}
Consider a baseleaf preserving homotopy $H:S^{1}_\Q\times [0,1]\rightarrow S^{1}_\Q$. By Lemma \ref{LimitPer}, on the baseleaf the homotopy has the following expression:
$$H_0(x,t)= q(t)x+h(x,t)$$
such that $t\mapsto q(t)\in\Q$ is a continuous function hence constant for $\Q$ is totally disconnected and we have the result.
\end{proof}

\bigskip

The following proposition gives the converse of Lemma \ref{LimitPer}.

\begin{prop}\label{LimPerExtMap}
For every rational $q$ and continuous limit periodic real valued function $h$, there is a unique continuous baseleaf preserving map $f$ from the solenoid to itself such that the following diagram commutes:

$$\xymatrix{S^{1}_\Q \ar[r]^{f} & S^{1}_\Q \\
			\R \ar[r]^{f_{0}} \ar@{^{(}->}[u]_{\nu} & \R \ar@{^{(}->}[u]_{\nu}}$$
where $f_{0}(x)= qx + h(x)$.
\end{prop}
\begin{proof}
Consider a rational $p/q$ such that $p$ and $q$ are coprime natural numbers and define $F:q\Z\times\R\rightarrow p\Z\times\R$ such that $F(qn, x)= (pn, f_{n}(x))$ where:
\begin{equation}\label{ext_f0}
f_{n}(x)= f_{0}(x+2\pi qn)-2\pi pn= \frac{p}{q}x + h(x+2\pi qn)
\end{equation}
for every integer $n$. Because $h$ is continuous and limit periodic, the function $g:q\Z\times\R\rightarrow \R$ such that $g(n,z)= h(x+2\pi n)$ admits a unique continuous extension $\hat{g}:q\hat{\Z}\times \R\rightarrow \R$ such that $\hat{g}(a,x+2\pi q)= \hat{g}(a+q,x)$. Then, there is a unique continuous extension $\hat{F}:q\hat{\Z}\times\R\rightarrow p\hat{\Z}\times\R$ of $F$ such that:
$$\hat{F}(qa,x)= \left(pa, \frac{p}{q}x + \hat{g}(qa, x)\right)$$
and satisfies the same structural condition as $F$:
$$\hat{F}(qa, x+2\pi q)= \hat{F}(q(a+1),x) + (-p,2\pi p)$$
By Lemma \ref{Exact_seq_II}, there is a continuous map $f$ such that the following diagram commutes:

$$\xymatrix{	S^{1}_{\Q} \ar[rrr]^{f} & & & S^{1}_{\Q} \\
			q\hat{\Z}\times\R \ar[rrr]^{\hat{F}} \ar[u]_{\exp} & & & p\hat{\Z}\times\R \ar[u]_{\exp} \\			
			\R \ar[rrr]^{(p/q)x+h(x)} \ar@{^{(}->}[u] & & & \R \ar@{^{(}->}[u]}$$
and we have the result.
\end{proof}

\begin{prop}
$\deg(f)= \deg(f')$ iff $f$ is homotopic to $f'$.
\end{prop}
\begin{proof}
If $f$ is homotopic to $f'$, then they have the same degree by Proposition \ref{Hom_Invariant}. Conversely, consider their expressions on the baseleaf: $f_0(x)= qx+h(x)$ and $f'_0(x)=qx+h'(x)$. Define the homotopy $H_0(x,t)= qx+th(x)+(1-t)h'(x)$. The construction given in Lemma \ref{LimPerExtMap} extends to a continuous family of maps hence there is a homotopy between $f$ and $f'$ whose expression in the baseleaf is $H_0$.
\end{proof}

\bigskip

Define the baseleaf preserving continuous map $z^{q}$ from the solenoid to itself as the lifting of the map $x\mapsto qx$ through the baseleaf.
\begin{obs}\label{Remark} Although they have the same notation, do not confuse this map with the one defined in the section \ref{adelic_sphere}: One is from the solenoid to itself and the other is from the (algebraic) solenoid to $\C^{*}$. The context will make the distinction clear. Among these maps, only the identity is leaf preserving and all the others permute leaves.\end{obs}

\begin{cor}
The homotopic classes of baseleaf preserving continuous maps from the solenoid to itself are in one to one correspondence with the rational numbers:
$$[S^{1}_\Q, S^{1}_\Q]_{L}\cong \Q$$
and the correspondence is given by the degree.
\end{cor}

\section{Adelic Riemann sphere}\label{adelic_sphere}

\begin{lema}\label{homotopy}
Consider a finite dimensional Lie group $G$ and a continuous function $f$ from the solenoid to $G$. There is $\varepsilon'>0$ such that for every $\varepsilon'>\varepsilon>0$, there is a continuous function $f'$ homotopic to $f$ such that $d(f,f')<\varepsilon$ and it factors through some canonical projection $\pi_n$.
\end{lema} 
\begin{proof}
Consider a bi--invariant metric on the Lie group and a convex neighborhood $V$ of the origin in the Lie algebra $\mathfrak{g}$ such that the exponential restricted to $V$ is a local homeomorphism. Consider $\varepsilon'>0$ such that $B_{\varepsilon'}(e)\subset \exp(V)$ and $\varepsilon>0$ such that $\varepsilon'>\varepsilon>0$. By Lemma \ref{Aprox}, there is a continuous function $f'$ such that $d(f,f')<\varepsilon$ and it factors through some canonical projection $\pi_n$.
Define the continuous function $a$ from the solenoid to $V\subset \mathfrak{g}$ such that $f'(p)f(p)^{-1}= \exp\left(a(p)\right)$. Define the homotopy:
$$H(p,t):= \exp\left(t\ a(p)\right)f(p)$$
Then $f'$ is homotopic to $f$ and we have the result.
\end{proof}
\bigskip

In the previous Lemma, the finite dimension hypothesis on the Lie group is too restrictive for we only need the existence of an exponential map. For example, in the infinite dimensional case, a loop group has exponential map but the Virasoro group does not.

If instead of the circle we take the inverse limit of the coverings of $\C^{*}$ just as before, we get the \textsf{algebraic solenoid} $\C^{*}_\Q$. Again, taking the inverse limit of the ramified coverings of the Riemann sphere $\C P^{1}$, we get the \textsf{adelic Riemann sphere} $\C P^{1}_\Q$. All of the ramifications occur at $0$ and $\infty$ hence topologically the adelic Riemann sphere is the topological suspension of the adelic solenoid.

The same theory in section \ref{adelic_solenoid} can be constructed for the algebraic solenoid almost verbatim. The complex structure of the algebraic solenoid as a laminated object is induced by the exponential map $\exp:\hat{\Z}\times \C\rightarrow \C^{*}_\Q$ with the obvious complex structure on $\hat{\Z}\times \C$.

A rank $n$ holomorphic vector bundle over the adelic Riemann sphere is a vector bundle whose clutching map is a holomorphic function from the algebraic solenoid to $GL(n,\C)$. We will denote simply as $z^{m/n}$ the map $\left(z\mapsto \pi_n(z)^{m}\right)$ from the algebraic solenoid to $\C^{*}$, where $m$ and $n$ are natural numbers (See remark \ref{Remark}). We will denote as $\I(q)$ the holomorphic line bundle whose clutching map is $z^{q}$ and $q$ is a rational number.

The following is the adelic version of the celebrated Birkhoff-Grothendieck Theorem:

\begin{teo}\label{BG}
Consider a rank $n$ holomorphic vector bundle $V$ over the adelic Riemann sphere $\C P^{1}_\Q$. There is a set of rational numbers $q_1,\ldots q_n$ such that:
$$V\cong \I(q_1)\oplus\ldots \I(q_n)$$
isomorphic as holomorphic vector bundles.
\end{teo}
\begin{proof}
The clutching map of $V$ is a holomorphic map $f$ from the algebraic solenoid $\C^{*}_\Q$ to $G= GL(n, \C)$ and consider its restriction $f'$ to the solenoid $S^{1}_\Q$. By Lemma \ref{homotopy}, there is a continuous homotopic map $f''\simeq f'$ that factors trough some canonical projection $\pi_m$; i.e. There is a continuous map $g: S^{1}\rightarrow G$ such that $f''= g\circ\pi_m$. Apply the $\boldsymbol{\pi}_1$ functor to this map:
$$\Z\cong \boldsymbol{\pi}_1(S^{1})\xrightarrow{\boldsymbol{\pi}_1(g)} \boldsymbol{\pi}_1(G)\cong \Z$$
Then, $\boldsymbol{\pi}_1(g)$ acts by multiplying by some integer $n$. Define the holomorphic map $\hat{g}: \C^{*}\rightarrow G$ such that $\hat{g}(z)= diag(z^{n},1,\ldots 1)$ and the pullback by the projection $f_m:= \hat{g}\circ \pi_m$. Define the holomorphic map $h:= f_m \cdot f^{-1}$. By construction, $h$ factors through the universal cover group $\tilde{G}\xrightarrow{\mu} G$ and the pullback of the lifting by the exponential map is a holomorphic map $a$ from the algebraic solenoid to the Lie algebra; i.e. $h(p)= \mu\left(\exp(a(p))\right)$ and $a(p)\in \mathfrak{g}$. We summarize the construction in the following commutative diagram:

$$\xymatrix{\mathfrak{g} \ar[rr]^{\exp} & & \tilde{G} \ar[d]^{\mu} \\
\C^{*}_\Q \ar[u]^{a} \ar[urr]^{\tilde{h}} \ar[rr]^{h} & & G }$$

\noi Define the holomorphic isotopy (continuous map such that it is holomorphic on each $t$):
$$H(p,t)= \mu\left(\exp(t\ a(p))\right)f(p)$$
Then, $f_m$ is holomorphic isotopic to $f$ and we have:
$$V\cong \I(n/m)\oplus\I(1/m)\ldots \I(1/m)\eqno\qedhere$$
\end{proof}

%\bigskip

In particular, for every holomorphic line bundle $L$ there is a rational number $q$ such that $L\cong \I(q)$ and, as the next Lemma shows, this number is unique.

\begin{lema}
Consider rational numbers $q$ and $q'$. Then, $\I(q)\cong \I(q')$ iff $q=q'$.
\end{lema}
\begin{proof}
There is a holomorphic isotopy $h$ between the clutching maps of $\I(q)$ and $\I(q')$; i.e. $h(p,0)= p^{q}$ and $h(p,0)= p^{q'}$ Following the proof of Theorem \ref{BG}, for every $t$ there is a holomorphic isotopy from $\left(p\mapsto h(p,t)\right)$ to $z^{q(t)}$ and these isotopies vary continuously with respect to $t$:
$$H(p, t, s):= \mu\left(\exp(s\ a(p,t))\right)h(p,t)$$
Then, $H(p,t,1)= p^{q(t)}$ is a holomorphic isotopy hence $(t\mapsto q(t))$ is a continuous map such that $q=q(0)$ and $q'=q(1)$. Because $\Q$ is totally disconnected, the map is constant and $q=q'$.
\end{proof}

\begin{defi}
For every holomorphic line bundle $L$ over the adelic Riemann sphere, define its \textsf{Chern character} $ch(L)$ as the unique rational number $q$ such that $L\cong \I(q)$.
\end{defi}

\begin{cor}
The Picard group of the adelic Riemann sphere is isomorphic to the rational additive abelian group:
$$Pic\left(\C P^{1}_\Q\right)\cong \Q$$
and the isomorphism is given by the Chern character.
\end{cor}
\begin{proof}
If $L\cong \I(q)$ and $L'\cong \I(q')$, then $L\otimes L'\cong \I(q)\otimes\I(q')\cong \I(q+q')$.
\end{proof}

\bigskip

The analog of the Picard group for higher rank vector bundles is the $K$--ring:

\begin{cor}
The $K$--ring of the adelic Riemann sphere is the following:
$$K\left(\C P^{1}_\Q\right)\cong \frac{\Z[x^{q}\ /\ q\in\Q]}{\left\langle\ (x^{q}-1)(x^{p}-1)\ /\ q,p\in\Q \right\rangle}$$
\end{cor}
\begin{proof}
The proof is almost verbatim to the usual one for $\C P^{1}$. By Theorem \ref{BG}, we have an ring epimorphism
$$\xi:\Z[z^{q}\ /\ q\in\Q]\rightarrow K\left(\C P^{1}_\Q\right)$$
such that $z^{q}\mapsto \I(q)$. Consider the homotopy:
$$F(z,t)= \left(\begin{array}{cc}
                           z^{a} & 0 \\
                           0 & 1  \\
                          \end{array}
                       \right)\left(\begin{array}{cc}
                           cos(\pi\ t/2) & -sen(\pi\ t/2) \\
                           sen(\pi\ t/2) & cos(\pi\ t/2)  \\
                          \end{array}
                       \right)\left(\begin{array}{cc}
                           z^{b} & 0 \\
                           0 & 1  \\
                          \end{array}
                       \right)\left(\begin{array}{cc}
                           cos(\pi\ t/2) & sen(\pi\ t/2) \\
                           -sen(\pi\ t/2) & cos(\pi\ t/2)  \\
                          \end{array}
                       \right)$$
Evaluating $t=0$ and $t=1$ we have that the following clutching functions are holomorphic isotopic:
$$\left(\begin{array}{cc}
                           z^{a+b} & 0 \\
                           0 & 1  \\
                          \end{array}
                       \right)\sim \left(\begin{array}{cc}
                           z^{a} & 0 \\
                           0 & z^{b}  \\
                          \end{array}
                       \right)$$
In particular we have the following isomorphism as holomorphic vector bundles:
$$\mathcal{O}(a+b)\oplus\mathcal{O}(0)\simeq \mathcal{O}(a)\oplus\mathcal{O}(b)$$
This implies that $(z^{a}-1)(z^{b}-1)= z^{a+b}-z^{a}-z^{b}+1$ is annihilated by the map $\xi$ for every pair $a,b$ of rational numbers and the epimorphism $\xi$ factors through the quotient:
$$\hat{\xi}: \frac{\Z[z^{q}\ /\ q\in\Q]}{\left\langle\ (z^{q}-1)(z^{p}-1)\ /\ q,p\in\Q \right\rangle}\rightarrow K\left(\C P^{1}_\Q\right)$$
To see that $\hat{\xi}$ is in fact an isomporhism, we just need to show an inverse. Define the map:
$$\nu: K\left(\C P^{1}_\Q\right)\rightarrow \frac{\Z[z^{q}\ /\ q\in\Q]}{\left\langle\ (z^{q}-1)(z^{p}-1)\ /\ q,p\in\Q \right\rangle}$$
as the ring morphism such that $\nu(\mathcal{O}(q))= z^{q}$. It is clear that $\nu$ is well defined and it is the inverse of $\hat{\xi}$.
\end{proof}

As it was expected from the continuous surjective map $\C P^{1}_\Q\rightarrow \C P^{1}$ and the fact that $K$ is a contravariant functor, we have the inclusion:
$$\frac{\Z[x]}{\left\langle\ (x-1)^{2} \right\rangle}\cong K(\C P^{1}) \hookrightarrow K\left(\C P^{1}_\Q\right)$$
where $x$ is the image of the class represented by the tautological complex line bundle. It is interesting to note the rational factors of the tautological class in the adelic case.

We say that a complex function of the adelic Riemann sphere is holomorphic if it is holomorphic on the algebraic solenoid and admits a continuous extension to the cusps.

\begin{prop}
A holomorphic function of the adelic Riemann sphere is constant:
$$Hol(\C P^{1}_\Q)\cong \C$$
\end{prop}
\begin{proof}
Consider a holomorphic function $f$ of the adelic Riemann sphere and its pullback $g$ by the baseleaf $\nu$. Then, $g$ is a holomorphic limit periodic respect to $x$ ($z=x+iy$) function such that:
$$\lim_{y\to +\infty} ||f(0)-g|_{Im(z)\geq y}\ ||_{\infty}=0$$
$$\lim_{y\to -\infty} ||f(\infty)-g|_{Im(z)\leq y}\ ||_{\infty}=0$$
Then, $g$ is a bounded entire function hence it is constant. Because the baseleaf is dense in the adelic Riemann sphere, by continuity $f$ is constant.
\end{proof}


\begin{thebibliography}{Coh93}

\bibitem[Be]{Bell}
H.Bell, K.R.Meyer, \emph{Limit Periodic Functions, Adding Machines, and Solenoids}, Journal of Dynamics and Differential Equations, vol.7, 3 (1995).

\bibitem[Bi]{Birkhoff}
G.D.Birkhoff, \emph{A theorem on matrices of analytic functions}, Math. Ann., vol.74, (1913), 122-133.


\bibitem[Bo]{Bohr}
H.Bohr, \emph{Zur Theorie der fastperiodischen Funcktionen I, II, III}, Acta Math., vol.45, (1924), 29-127.

\bibitem[Fe]{Feigenbaum}
M.Feigenbaum, \emph{Universal Behavior in Nonlinear Systems}, Los Alamos science, vol.1, (1980), 4-27,

\bibitem[Gr]{Grothendieck}
A.Grothendieck, \emph{Sur la classification des fibres holomorphes sur la sph\`ere de
Riemann}, Amer. J. Math., vol.79, (1957), 121-138.

\bibitem[SGA1]{SGA1}
A.Grothendieck, M.Raynaud, \emph{Rev\^etement Etales et Groupe Fondamental (SGA1)}, Lecture Note in
Math., Springer, Berlin Heidelberg New York, vol.224, (1971).


\bibitem[LM]{LM}
M. Lyubich, Y. Minsky, \emph{Laminations in holomorphic dynamics},  J. Differential Geom., vol.47, 1, (1997).


\bibitem[Od]{Odden}
C.Odden, \emph{The baseleaf preserving mapping class group of the universal hyperbolic solenoid}, Trans. Amer. Math. Soc., vol.357, (2005), 1829-1858.


\bibitem [RV]{RV}
D.Ramakrishnan, R.Valenza, \emph{Fourier Analysis on Number Fields}.  Graduate Texts in Mathematics {\bf 186}, Springer-Verlag,
New York, (1999).

\bibitem[Su]{Universalities}
D.Sullivan, \emph{Linking the universalities of Milnor-Thurston, Feigenbaum and Ahlfors-Bers}, L. R. Goldberg and A. V. Phillips, editors, Topological Methods in Modern Mathematics, Publish or Perish, (1993), 543-563.

\bibitem [Su2]{Su2} 
D.Sullivan, \emph{Bounds, quadratic differentials, and renormalization conjectures}, American Mathematical Society centennial publications, (Providence, RI, 1988), Amer. Math. Soc., Providence, RI, vol.2, (1992), 417-466.

\bibitem [Ta] {Ta}
J.T.Tate, \emph{Fourier analysis in number fields, and Hecke's zeta-functions}, Algebraic Number Theory, ed. Cassels, J.W.S. \& Fr\"{o}hlich, A., Cambridge U. Press, Cambridge UK, (2010), 305-347.

\end{thebibliography}
\end{document}